\documentclass[11pt]{article}
\usepackage{latexsym,amssymb,amsmath,a4wide,theorem}
\usepackage[isolatin]{inputenc}
\usepackage[T1]{fontenc}
\usepackage{times}
\usepackage{cite}
\usepackage{amsfonts}
\usepackage{hyperref}

\def\sqr#1#2{{\vbox{\hrule height.#2pt
     \hbox{\vrule width.#2pt height#1pt \kern#1pt
           \vrule width.#2pt}
     \hrule height.#2pt}}}

\newtheorem{theorem}{Theorem}[section]
\newtheorem{lemma}[theorem]{Lemma}

\title{\Large\bf On the existence of multiple normalized solutions for a class of fractional Choquard equations with mixed nonlinearities
}
\author{
{Yongpeng Chen}\thanks{Email:yongpengchen@mail.bnu.edu.cn}\\
\small School of Science, Guangxi University of Science and Technology, Liuzhou, 545006. P.R.China.\\
{Zhipeng Yang}\thanks{Email:yangzhipeng326@163.com}\\
\small Department\, of \,Mathematics, Yunnan\, Normal\, University, Kunming, 650500, P.R.China.\\
{Jianjun Zhang}\thanks{Corresponding author:zhangjianjun09@tsinghua.org.cn}\\
\small College of Mathematics and Statistics, Chongqing Jiaotong University, Chongqing 400074, PR China\\
}
\date{}

\begin{document}
\maketitle
\begin{abstract}
We investigate the existence of normalized solutions for the following nonlinear fractional Choquard equation:
$$
(-\Delta)^s u+V(\epsilon x)u=\lambda u+\left(I_\alpha *|u|^q\right)|u|^{q-2} u+\left(I_\alpha *|u|^p\right)|u|^{p-2} u, \quad x \in \mathbb{R}^N,
$$
subject to the constraint
$$
\int_{\mathbb{R}^N}|u|^2 \mathrm{d}x=a>0,
$$
where $N>2 s, s \in(0,1), \alpha \in(0, N), \frac{N+\alpha}{N}<q<\frac{N+2 s+\alpha}{N}<p\leq \frac{N+\alpha}{N-2 s}$,  $\epsilon>0$ is a parameter, and $\lambda \in \mathbb{R}$ serves as an unknown parameter acting as a Lagrange multiplier. By employing the Lusternik-Schnirelmann category theory, we estimate the number of normalized solutions to this problem by virtue of the category of the set of minimum points of the potential function $V$.
\end{abstract}


{\bf AMS} Subject Classification: 35A15, 35B40, 35J20
\numberwithin{equation}{section}
\section{Introduction and main results}
In this paper, we aim to explore  the existence of multiple normalized solutions for the following nonlinear fractional Choquard equation:
\begin{equation}\label{eq1}
\begin{cases}(-\Delta)^s u+V(\epsilon x)u=\lambda u+\left(I_\alpha *|u|^q\right)|u|^{q-2} u+\left(I_\alpha *|u|^p\right)|u|^{p-2} u, & x \in \mathbb{R}^N, \\
\int_{\mathbb{R}^N}|u|^2  \mathrm{d}x=a.
\end{cases}
\end{equation}
Here $N>2 s, s \in(0,1), \alpha \in(0, N)$, and $\frac{N+\alpha}{N}<q<\frac{N+2 s+\alpha}{N}<p\leq 2_{\alpha, s}^*$, where $2_{\alpha, s}^*:=\frac{N+\alpha}{N-2 s}$ is the upper Hardy-Littlewood-Sobolev critical exponent.  Furthermore,  $\epsilon>0$ is a parameter, and $\lambda \in \mathbb{R}$  is an unknown parameter acting as a Lagrange multiplier. Additionally, the Riesz potential $I_\alpha$ is defined for each point $x\in\mathbb{R}^N \backslash\{0\}$ by the formula
$$
I_\alpha(x):=\frac{A_\alpha}{|x|^{N-\alpha}}, \quad\text {where } A_\alpha=\frac{\Gamma\left(\frac{N-\alpha}{2}\right)}{\Gamma\left(\frac{\alpha}{2}\right) \pi^{\frac{N}{2}} 2^\alpha}
$$
with $\Gamma$ representing the Gamma function. The fractional Laplace operator $(-\Delta)^s$ in $\mathbb{R}^N$  is defined for functions in the Schwartz class as follows:
$$
(-\Delta)^s u(x)=\mathcal{F}^{-1}\left(|\xi|^{2 s} \mathcal{F}(u)\right)(x),\quad x \in \mathbb{R}^N,
$$
where $\mathcal{F}$ and $\mathcal{F}^{-1}$ represent the Fourier transform and its inverse, respectively. Alternatively, the fractional Laplacian can also be expressed in the following two forms:
\begin{equation*}
\begin{aligned}
(-\Delta)^s u(x)=&C(N, s) \mathrm{P.V.} \int_{\mathbb{R}^N} \frac{u(x)-u(y)}{|x-y|^{N+2 s}} \mathrm{d} y\\
=&-\frac{C(N, s)}{2} \int_{\mathbb{R}^N} \frac{u(x+y)-u(x-y)-2 u(x)}{|y|^{N+2 s}} \mathrm{d} y, \quad x \in \mathbb{R}^N,
\end{aligned}
\end{equation*}
where $\mathrm{P.V.}$ denotes the Cauchy principal value of the integral, and $C(N, s)$ is a normalization constant (see \cite{MR2944369} for details).
The fractional Laplacian was first introduced by Laskin in \cite{MR1755089} as an extension of the Feynman path integral, transitioning from Brownian quantum-mechanical paths to those resembling L\'evy processes. It plays a critical role in describing phenomena like Bonson stars and modeling water wave dynamics \cite{MR2318846}. Moreover, this operator has significant applications across diverse fields, including biology, chemistry, and finance, as highlighted in works like \cite{MR3009717, MR2042661, MR2737789, MR1809268, MR2270163}, and has been the subject of extensive research over the past few decades \cite{fra-1, fra-2, fra-3, fra-4, fra-5, fra-6}.

In the study of normalized solutions for the fractional Choquard equation:
$$
(-\Delta)^s u+V(x)u=\lambda u+\left(I_\alpha *|u|^p\right)|u|^{p-2} u, \quad x \in \ \mathbb{R}^N,
$$
three important quantitative indices play a crucial role: the Hardy-Littlewood-Sobolev lower critical exponent $\frac{N+\alpha}{N}$, the $L^2$-critical exponent $\frac{N+2 s+\alpha}{N}$, and the Hardy-Littlewood-Sobolev upper critical exponent $2_{\alpha, s}^*=\frac{N+\alpha}{N-2 s}$. The range of $p$ in relation to these indices defines distinct research approaches and presents various challenges.

Cingolani et al. \cite{Cingolani} investigated the existence of solutions to the $L^2$-
subcritical nonlocal  problem:
\begin{equation*}
\begin{cases}
(-\Delta)^s u+\mu u=\left(I_\alpha * F(u)\right) f(u) \quad \text {in}\quad \mathbb{R}^N, \\
\int_{\mathbb{R}^N} u^2 \mathrm{d}x=c,
\end{cases}
\end{equation*}
where $s \in(0,1), N \geq 2, \alpha \in(0, N), F \in C^1(\mathbb{R}, \mathbb{R})$ with $f=F^{\prime}, c>0$, and $\mu$ is a Lagrange multiplier. If $f$ satisfies the conditions of almost optimal $L^2$-subcritical growth, then, through the application of new deformation arguments, there exists a value $c_0 \geq 0$ such that for any $c>c_0$, the problem  possesses a radially symmetric solution.

Li and Luo \cite{li-luo} examined the following nonlinear $L^2$-
supcritical fractional Choquard equations:
\begin{equation*}
\begin{cases}
(-\Delta)^s u-\lambda u=\left(I_\alpha *|u|^p\right)|u|^{p-2} u\quad  \text {in }\  \mathbb{R}^{N},\\
\int_{\mathbb{R}^N} u^2 \mathrm{d} x=c,
\end{cases}
\end{equation*}
where $N \geq 3, s \in(0,1), \alpha \in(0, N), p \in\left(\max \left\{1+\frac{\alpha+2 s}{N}, 2\right\}, \frac{N+\alpha}{N-2 s}\right)$. By applying the constrained minimization approach, they established that this problem has a solution with minimal energy for any positive value of $c>0$. They also provided insights into the limiting behavior of these constrained solutions as $c$ tends toward both zero and infinity. Additionally, using a minimax technique, it was shown that there are infinitely many radial solutions available for every $c>0$.

Yu et al. \cite{yu-tang-zhang} studied  the following lower critical fractional Choquard equation:
\begin{equation*}
\begin{cases}
(-\Delta)^s u=\lambda u+\gamma\left(I_\alpha *|u|^{1+\frac{\alpha}{N}}\right)|u|^{\frac{\alpha}{N}-1} u+\mu|u|^{q-2} u \text { in } \mathbb{R}^N,\\
\int_{\mathbb{R}^N}|u|^2 \mathrm{d}x=a^2,
\end{cases}
\end{equation*}
where $N \geq 3, s \in(0,1), \alpha \in(0, N), a, \gamma, \mu>0$ and $2<q \leq 2_s^*:=2 N /(N-2 s)$.
When $2<q<2+4s/N$, the authors demonstrated the existence of radially normalized ground states by employing an extremal function construction technique. As for $2+4s/N<q<2_s^*$, they proposed a homotopy-stable family to show both the existence of Palais-Smale sequences and their compactness, leading to the conclusion regarding normalized ground states. Finally, for $q=2_s^*$, they established the presence of normalized ground states through Sobolev subcritical approximation methods.

He et al. \cite{he-1} examined the upper critical fractional Choquard equation:
\begin{equation*}
\begin{cases}
(-\Delta)^s u=\lambda u+\mu|u|^{q-2} u+\left(I_\alpha *|u|^{2_{\alpha, s}^*}\right)|u|^{2_{\alpha, s}^*-2} u, \quad x \in \mathbb{R}^N,\\
\int_{\mathbb{R}^N} u^2 \mathrm{~d} x=a^2,
\end{cases}
\end{equation*}
where $s \in(0,1), N>$ $2 s, 0<\alpha<\min \{N, 4 s\}, 2<q<2_s^*$, and $\mu \in \mathbb{R}$. The authors provided several results regarding existence and non-existence of solutions while also exploring the qualitative behavior of ground states as  $\mu \rightarrow 0^{+}$.

For additional insights into the normalized solutions of fractional Choquard equations, we suggest referring to \cite{feng-1, tao, feng-2, chen, liu, guo, lan, MR4767698,MR4391690} and their associated references. A review of the current literature on normalized solutions of fractional Choquard equations reveals that there are few studies leveraging the characteristics of the potential function to clarify the number of solutions. In contrast, significant progress has been made in examining normalized solutions for Schr\"{o}dinger equations. For instance, Alves et al. \cite{alves-thin} explored the existence of multiple normalized solutions within a specific class of elliptic problems represented by the following system:
$$
\begin{cases}
-\Delta u+V(\epsilon x) u=\lambda u+f(u), \quad x \in \mathbb{R}^N,\\ \int_{\mathbb{R}^N}|u|^2 \mathrm{d} x=a^2,
\end{cases}
$$
where the potential function $V$ meets the following criteria:
\begin{itemize}
\item[$(V_1)$] $V \in C\left(\mathbb{R}^N, \mathbb{R}\right) \cap L^{\infty}\left(\mathbb{R}^N\right)$ and $V(0)=0$,
\item[$(V_2)$] $0=\inf\limits_{x \in \mathbb{R}^N} V(x)<\liminf\limits_{|x| \rightarrow+\infty} V(x)=V_{\infty}$,
\end{itemize}
and $f$ is a differentiable function of $L^2$-subcritical growth. Utilizing minimization methods along with the Lusternik-Schnirelmann category, they demonstrated that the quantity of normalized solutions is associated with the topology of the set where the potential function $V$ attains its minimum value.

Li et al. \cite{li-xu-zhu} considered the following Schr\"{o}dinger equation:
$$
\left\{\begin{array}{l}
-\Delta u=\lambda u+h(\epsilon x)|u|^{q-2} u+\eta|u|^{p-2} u, \quad x \in \mathbb{R}^N, \\
\int_{\mathbb{R}^N}|u|^2 \mathrm{~d} x=a^2,
\end{array}\right.
$$
where $a, \epsilon, \eta>0$, $q$ is $L^2$-subcritical, and $p$ is $L^2$-supercritical. Here $h$ is a positive and continuous function that satisfies
\[ h_{\infty}:=\lim _{|x| \rightarrow+\infty} h(x)<\max\limits_{x\in\mathbb{R}^N}h(x).\]
It has been established that when $\epsilon$ is sufficiently small, the number of normalized solutions is at least equal to the count of global maximum points of $h$. For additional findings on this topic, we direct readers to references such as  \cite{alves-4, alves-3, alves-2,Bartsch-2, Ikoma-2,ding-2}.

Based on the previous discussion and primarily inspired by the results in \cite{alves-thin},  \cite{li-xu-zhu} and \cite{meng-he}, we aim to tackle problem \eqref{eq1} and demonstrate the existence of multiple normalized  solutions through the application of the Lusternik-Schnirelmann category theory.

Before presenting the primary findings, we define the following functional spaces. For any $s \in(0,1)$, the Sobolev space $H^s\left(\mathbb{R}^N\right)$ is characterized as follows:
$$
\begin{aligned}
H^s\left(\mathbb{R}^N\right) & =\left\{u \in L^2\left(\mathbb{R}^N\right): \frac{|u(x)-u(y)|}{|x-y|^{\frac{N}{2}+s}} \in L^2\left(\mathbb{R}^N \times \mathbb{R}^N\right)\right\} \\
& =\left\{u \in L^2\left(\mathbb{R}^N\right): \int_{\mathbb{R}^N}\left(1+|\xi|^{2 s}\right)|\mathcal{F}(u)(\xi)|^2<\infty\right\}
\end{aligned}
$$
with its norm defined by
$$
\left(\int_{\mathbb{R}^N} \int_{\mathbb{R}^N} \frac{|u(x)-u(y)|^2}{|x-y|^{N+2 s}} \mathrm{d} x \mathrm{d}y+\int_{\mathbb{R}^N}u^2\mathrm{d}x\right)^{\frac{1}{2}}.
$$
For $u \in H^s\left(\mathbb{R}^N\right)$, applying Propositions 3.4 and 3.6 from \cite{MR2944369}, we can derive
$$
\int_{\mathbb{R}^N}|(-\Delta)^{\frac{s}{2}} u|^2\mathrm{d}x=\int_{\mathbb{R}^N}|\xi|^{2 s}|\mathcal{F}(u)(\xi)|^2 \mathrm{d} \xi=\frac{1}{2} C(N, s) \int_{\mathbb{R}^N} \int_{\mathbb{R}^N}\frac{|u(x)-u(y)|^2}{|x-y|^{N+2 s}} \mathrm{d} x \mathrm{d} y.
$$
Thus, we will utilize the norm denoted as
$$
\|u\|_{H^s\left(\mathbb{R}^N\right)}:=\left(\int_{\mathbb{R}^N}|u|^2\mathrm{d}x+\int_{\mathbb{R}^N}|(-\Delta)^{\frac{s}{2}} u|^2\mathrm{d}x\right)^{\frac{1}{2}},
$$
and its associated inner product represented by
$$
\left\langle u, v\right\rangle:=\int_{\mathbb{R}^N}(-\Delta)^{\frac{s}{2}} u(-\Delta)^{\frac{s}{2}} v \mathrm{d} x+\int_{\mathbb{R}^N} uv\mathrm{d}x .
$$
Additionally, we introduce the space
\[\mathcal{D}^{s, 2}\left(\mathbb{R}^N\right)=\left\{u \in L^{2_s^*}\left(\mathbb{R}^N\right): \int_{\mathbb{R}^N} \int_{\mathbb{R}^N} \frac{|u(x)-u(y)|^2}{|x-y|^{N+2 s}} \mathrm{d} x \mathrm{d} y<+\infty\right\}
\]
equipped with the norm
\[
\|u\|:=\left(\int_{\mathbb{R}^N}|(-\Delta)^{\frac{s}{2}} u|^2\mathrm{d}x\right)^{\frac{1}{2}}.
\]

To obtain the normalized solutions for equation \eqref{eq1}, it is essential to examine the critical points of the energy functional defined as follows:
\begin{equation}\label{energyfunctional}
\begin{aligned}
J_\epsilon(u)=&\frac{1}{2}\int_{\mathbb{R}^N}|(-\Delta)^{\frac{s}{2}} u|^2\mathrm{d}x+\frac{1}{2}\int_{\mathbb{R}^N}V(\epsilon x)u^2\mathrm{d}x-\frac{1}{2q} \int_{\mathbb{R}^N}\left(I_\alpha *|u|^q\right)|u|^q \mathrm{d} x\\
&-\frac{1}{2p}\int_{\mathbb{R}^N}\left(I_\alpha *|u|^p\right)|u|^p\mathrm{d} x,
\end{aligned}
\end{equation}
subject to the constraint
$$
S(a)=\left\{u \in H^s\left(\mathbb{R}^N\right): \int_{\mathbb{R}^N}|u|^2 \mathrm{d}x=a\right\}.
$$

\begin{lemma}\cite{lieb-loss} \label{le1.1} Let $r, t>1$ and $\alpha \in(0, N)$ with $\frac{1}{r}+\frac{1}{t}=1+\frac{\alpha}{N}$. Let $f \in L^r\left(\mathbb{R}^N\right)$ and $h \in L^t\left(\mathbb{R}^N\right)$. Then there exists a sharp constant $C(r, t, \alpha, N)$ independent of $f, h$ satisfying
$$
\int_{\mathbb{R}^N} \int_{\mathbb{R}^N} \frac{f(x) h(y)}{|x-y|^{N-\alpha}} \mathrm{d}x \mathrm{d} y \leq C(r, t, \alpha, N)\|f\|_r\|h\|_t .
$$

\end{lemma}

According to this lemma, the functional $J_\epsilon$ defined in \eqref{energyfunctional} is well-defined within the space $H^s\left(\mathbb{R}^N\right)$ and is of class $ C^1$. For convenience, we define
$$Q(u):=\int_{\mathbb{R}^N}\left(I_\alpha *|u|^q\right)|u|^q \mathrm{d} x$$  and $$P(u):=\int_{\mathbb{R}^N}\left(I_\alpha *|u|^p\right)|u|^p\mathrm{d} x.$$
Thus, the functional $J_\epsilon$ can be expressed as follows:
$$J_\epsilon(u)=\frac{1}{2}\|u\|^2+\frac{1}{2}\int_{\mathbb{R}^N}V(\epsilon x)u^2dx-\frac{1}{2q} Q(u)-\frac{1}{2p}P(u).$$

\begin{lemma}\cite{feng-zhang}\label{le1.2}  Let $N>2 s, 0<s<1$ and $\frac{N+\alpha}{N}<t<2_{\alpha, s}^*$, then, for all $u \in H^s\left(\mathbb{R}^N\right)$,
$$
\int_{\mathbb{R}^N}\left(I_\alpha *|u|^t\right)|u|^t \mathrm{d}x \leq C_{\alpha, t}\left\|u\right\|^{2 t \gamma_{t, s}}\|u\|_2^{2 t\left(1-\gamma_{t, s}\right)},
$$
where $\gamma_{t, s}=\frac{N t-N-\alpha}{2 t s}$ and the optimal constant $C_{\alpha, t}$ is given by
$$
C_{\alpha, t}=\frac{2 s t}{2 s t-N t+N+\alpha}\left(\frac{2 s t-N t+N+\alpha}{N t-N-\alpha}\right)^{\frac{Nt-N-\alpha}{2 s}}\|U\|_2^{2-2 t} \text {, }
$$
where $U$ is the ground state solution of
$$
(-\Delta)^s U+U-\left(I_\alpha *|U|^t\right)|U|^{t-2} U=0,\,x\in\mathbb{R}^N .
$$
\end{lemma}

\begin{lemma}\cite{he-r}\label{le1.3}
Let $\mathcal{S}_\alpha$ be the best constant
$$
\mathcal{S}_\alpha:=\inf _{u \in \mathcal{D}^{s, 2}\left(\mathbb{R}^N\right) \backslash\{0\}} \frac{\int_{\mathbb{R}^N}\left|(-\Delta)^{\frac{s}{2}} u\right|^2 \mathrm{d} x}{\left(\int_{\mathbb{R}^N}\left(I_\alpha *|u|^{2_{\alpha,s}^*}\right)|u|^{2_{\alpha,s}^*} \mathrm{d} x\right)^{\frac{1}{2_{\alpha,s}^*}}},
$$
and $\mathcal{S}_\alpha$ is achieved if and only if $u$ is of the form
$$
C\left(\frac{\varepsilon}{\varepsilon^2+\left|x-x_0\right|^2}\right)^{\frac{N-2 s}{2}}, \quad x \in \mathbb{R}^N,
$$
for some $x_0 \in \mathbb{R}^N, C>0$ and $\varepsilon>0$.
\end{lemma}

Let $\gamma_{2_{\alpha,s}^*, s}=1$ and define $C_{\alpha, 2_{\alpha,s}^*}=\left(\frac{1}{\mathcal{S}_\alpha}\right)^{2_{\alpha,s}^*}$. Then, based on Lemmas \ref{le1.2} and \ref{le1.3}, we get 
\begin{equation}\label{1.3}
\int_{\mathbb{R}^N}\left(I_\alpha *|u|^t\right)|u|^t d x \leq C_{\alpha, t}\left\| u\right\|^{2 t \gamma_{t, s}}\|u\|_2^{2 t\left(1-\gamma_{t, s}\right)},
\end{equation}
where $t\in\left(\frac{N+\alpha}{N},2_{\alpha,s}^*\right]$. Next, we introduce
$$
\mathcal{K}:=\frac{p \gamma_{p,s}-q \gamma_{q,s}}{1-q \gamma_{q,s}}\left(\frac{1-q \gamma_{q,s}}{p \gamma_{p,s}-1}\right)^{\frac{p \gamma_{p,s}-1}{p \gamma_{p,s}-q \gamma_{q,s}}}\left(\frac{C_{\alpha, q}}{q}\right)^{\frac{p \gamma_{p,s}-1}{p \gamma_{p,s}-q \gamma_{q,s}}}\left(\frac{C_{\alpha, p}}{p}\right)^{\frac{1-q \gamma_{q,s}}{p \gamma_{p,s}-q \gamma_{q,s}}},
$$
and impose the following conditions on $a$:
\begin{equation}\label{1.4}
a^{\frac{q(1-\gamma_{q,s})(p \gamma_{p,s}-1)+p(1-\gamma_{p,s})(1-q \gamma_{q,s})}{p \gamma_{p,s}-q \gamma_{q,s}}}<\frac{1}{\mathcal{K}},
\end{equation}
and
\begin{equation}\label{1.5}
\begin{aligned}
\left( a^{q\left(1-\gamma_{q,s}\right)}\right)^{\frac{1}{2(2^*_{\alpha,s}-q \gamma_q)}} \leq\left(\frac{\left(2-2q \gamma_{q,s}\right) C_{\alpha, q} 2_{\alpha,s}^* \mathcal{S}_\alpha^{2_{\alpha,s}^*}}{q(22^*_{\alpha,s}-2)}\right)^{\frac{-1}{2(2^*_{\alpha,s}-q \gamma_q)}} \mathcal{S}_\alpha^{\frac{N+\alpha}{2(\alpha+2s)}} .
\end{aligned}
\end{equation}
In this article, we assume that the condition \eqref{1.4} is satisfied. Before presenting the main theorem, let us clarify some notations.  For $\delta>0$, define
$$
\mathcal{M}=\left\{x \in \mathbb{R}^N: V(x)=0\right\}
$$
and
$$
\mathcal{M}_\delta=\left\{x \in \mathbb{R}^N: \operatorname{dist}(x, \mathcal{M}) \leq \delta\right\} .
$$
Now, we can present our primary findings in the following manner.
\begin{theorem}\label{thm1}
Assume $\left(V_1\right)$-$\left(V_2\right)$ and \eqref{1.4}, with the additional assumption that \eqref{1.5} holds if $p=2_{\alpha,s}^*$, then for each $\delta>0$, there exist constants $\epsilon_0>0$ and $V_*>0$ such that, for $0<\epsilon<\epsilon_0$ and $\|V\|_{\infty}<V_*$, the equation \eqref{eq1} has at least $cat_{\mathcal{M}_\delta}(\mathcal{M})$ pairs of weak solutions $\left(u_i, \lambda_i\right) \in H^s\left(\mathbb{R}^N\right) \times \mathbb{R}$ under the constraints of $\|u_i\|_2^2=a$. Moreover, $\lambda_i<0$, and $J_\epsilon(u_i)<0$.
\end{theorem}

\textbf{Notation.}~In this paper we make use of the following notations.
\begin{itemize}
\item[$\bullet$] $B_{R}(x)$ represents the open ball with a radius of $R$ centered at $x$, where $R>0$ and  $x\in\mathbb{R}^N$.
\item[$\bullet$] The symbols $C$ and $C_i$, where $i\in\mathbb{N}^+$, represent positive constants in the following context.
\item[$\bullet$]  "$\rightarrow$" and "$\rightharpoonup$" denote strong convergence and weak convergence, respectively.
\item[$\bullet$] $o_n(1)$ is a quantity tending to 0 as $n\to \infty$.
\item[$\bullet$] $\|u\|_r=(\int_{\mathbb{R}^N}|u|^r\mathrm{d}x)^{\frac{1}{r}}$ denotes the norm of $u$ in
$L^r(\mathbb{R}^N)$ for $r\geq 1$.
\item[$\bullet$] $\|u\|_\infty=\text{ess}\sup\limits_{x\in\mathbb{R}^N}|f(x)|$ denotes the norm of $u$ in $L^\infty(\mathbb{R}^N)$.

\item[$\bullet$] $2_s^*:=\frac{2 N}{N-2 s}$ is the fractional critical Sobolev exponent.


\end{itemize}

\section{Preliminaries}

For every $u \in S(a)$, according to \eqref{1.3}, we can see that
\begin{equation}\label{2.1}
\begin{aligned}
 J_\epsilon(u) &\geq \frac{1}{2}\|u\|^2-\frac{1}{2 q} C_{\alpha, q}\|u\|^{2q \gamma_{q, s}}\|u\|_2^{2 q\left(1-\gamma_{q, s}\right)}-\frac{1}{2 p} C_{\alpha, p}\|u\|^{2 p \gamma_{p, s}}\|u\|_2^{2 p\left(1-\gamma_{p, s}\right)} \\
& =\frac{1}{2}\|u\|^2\left(1-\frac{ a^{q\left(1-\gamma_{q,s}\right)}}{q} C_{\alpha, q}\|u\|^{2 q \gamma_{q, s}-2}-\frac{ a^{p\left(1-\gamma_{p, s}\right)}}{p} C_{\alpha, p}\|u\|^{2 p \gamma_{p, s}-2}\right) \\
& =g_a(\|u\|),
\end{aligned}
\end{equation}
where
$$
g_a(r):=\frac{1}{2}r^2-\frac{ a^{q\left(1-\gamma_{q,s}\right)}}{2q} C_{\alpha, q}r^{2 q \gamma_{q, s}}-\frac{ a^{p\left(1-\gamma_{p, s}\right)}}{2p} C_{\alpha, p}r^{2 p \gamma_{p, s}}, \quad r>0 .
$$
Define $g_a(r)=\frac{1}{2}r^2 w_a(r)$, where
\begin{equation}\label{2.2}
w_a(r):=1-\frac{ a^{q\left(1-\gamma_{q,s}\right)}}{q} C_{\alpha, q}r^{2 q \gamma_{q, s}-2}-\frac{ a^{p\left(1-\gamma_{p, s}\right)}}{p} C_{\alpha, p}r^{2 p \gamma_{p, s}-2}, \quad r>0 .
\end{equation}
For $t\in\left(\frac{N+\alpha}{N},2_{\alpha,s}^*\right]$, by the definition of $\gamma_{t,s}$, we have
$$
t \gamma_{t,s}\left\{\begin{array}{ll}
<1, & \text{if}\quad\frac{N+\alpha}{N}<t<\frac{N+\alpha+2s}{N}, \\
=1, & \text{if}\quad t=\frac{N+\alpha+2s}{N}, \\
>1, & \text{if}\quad\frac{N+\alpha+2s}{N}<t \leq 2_{\alpha,s}^*.
\end{array}\right.
$$
Therefore, from \eqref{2.2}, $\lim\limits _{r \rightarrow 0^{+}} w_a(r)=-\infty$ and $\lim\limits _{r \rightarrow+\infty} w_a(r)=-\infty$. Upon direct calculation,
$$
w_a^{\prime}(r)=-\frac{ a^{q\left(1-\gamma_{q,s}\right)}}{q} C_{\alpha, q}(2 q \gamma_{q, s}-2)r^{2 q \gamma_{q, s}-3}-\frac{ a^{p\left(1-\gamma_{p, s}\right)}}{p} C_{\alpha, p}(2 p \gamma_{p, s}-2)r^{2 p \gamma_{p, s}-3}.
$$
Then, the equation $w_a^{\prime}(r)=0$ possesses a unique solution given by:
$$
r_0=\left(\frac{\left(2-2q \gamma_{q,s}\right) \frac{1}{q}  C_{\alpha, q} a^{q\left(1-\gamma{q,s}\right)}}{\left(2p \gamma_{p,s}-2\right) \frac{1}{p} C_{\alpha, p} a^{p\left(1-\gamma_{p,s}\right)}}\right)^{\frac{1}{2(p \gamma_{p,s}-q \gamma_{q,s})}}.
$$
Furthermore, the maximum value of $w_a(r)$ is identified as follows:
$$
w_a\left(r_0\right)=1-\mathcal{K}a^{\frac{q(1-\gamma_{q,s})(p \gamma_{p,s}-1)+p(1-\gamma_{p,s})(1-q \gamma_{q,s})}{p \gamma_{p,s}-q \gamma_{q,s}}}.
$$
Thus, according to \eqref{1.4}, the maximum of $w_a(r)$ is positive and $w_a(r)=0$ has precisely two roots such that $0<R_0<R_1<\infty$. In conclusion, the properties of $w_a(r)$ can be summarized as follows:
$$
\left\{\begin{array}{l}
\lim\limits _{r \rightarrow 0^{+}} w_a(r)=-\infty \ \text{and}\  \lim\limits _{r \rightarrow+\infty} w_a(r)=-\infty;\\
w_a\left(R_0\right)=w_a\left(R_1\right)=0, \  \text {and}\  w_a\left(r_0\right)>0;  \\
w_a(r)\  \text{is increasing in }\left(0, r_0\right) \text { and   decreasing in} \left(r_0, +\infty\right).
\end{array}\right.
$$
Noting that $g_a(r)=\frac{1}{2}r^2 w_a(r)$, we can observe that   $g_a(r)$ possesses the following properties:
$$
\left\{\begin{array}{l}
g_a(0)=g_a\left(R_0\right)=g_a\left(R_1\right)=0 ;  \\
\lim\limits _{r \rightarrow+\infty} g_a(r)=-\infty ;  \\
g_a\left(r\right)<0 \ \text{for}\ r \in\left(0, R_0\right) \text { and } g_a\left(r\right)>0 \ \text{for}\  r \in\left(R_0, R_1\right).
\end{array}\right.
$$
Now, let $\tau:(0,+\infty) \rightarrow[0,1]$ be a non-increasing and $C^{\infty}$ function such that
$$
\tau(x)=\left\{\begin{array}{ll}
1, & \text { if } x \leq R_0, \\
0, & \text { if } x \geq R_1.
\end{array}\right.
$$
Consider the truncated functional:
$$
J_{\epsilon, T}(u):=  \frac{1}{2}\|u\|^2+\frac{1}{2}\int_{\mathbb{R}^N}V(\epsilon x)u^2\mathrm{d}x-\frac{1}{2q}Q(u) -\tau\left(\| u\|\right)\frac{1}{2p}P(u).
$$
It is easy to see that $J_{\epsilon, T} \in C^1\left(H^s\left(\mathbb{R}^N\right), \mathbb{R}\right)$. According to \eqref{1.3}, for any $u \in S(a)$,
$$
\begin{aligned}
J_{\epsilon, T}(u) &\geq \frac{1}{2}\|u\|^2-\frac{1}{2 q} C_{\alpha, q}\|u\|^{2q \gamma_{q, s}}a^{ q\left(1-\gamma_{q, s}\right)}-\frac{1}{2 p} C_{\alpha, p}\tau\left(\| u\|\right)\|u\|^{2 p \gamma_{p, s}}a^{ p\left(1-\gamma_{p, s}\right)} \\
&=g_{a,T}\left(\| u\|\right),
\end{aligned}
$$
where
$$
g_{a,T}(r):=\frac{1}{2}r^2-\frac{1}{2 q} C_{\alpha, q}r^{2q \gamma_{q, s}}a^{ q\left(1-\gamma_{q, s}\right)}-\frac{1}{2 p} C_{\alpha, p}\tau\left(r\right)r^{2 p \gamma_{p, s}}a^{ p\left(1-\gamma_{p, s}\right)}.
$$
By the properties of $g_a(r)$ and $\tau(r)$,  $g_{a,T}(r)$ possesses the following properties:
$$
\left\{\begin{array}{l}
g_{a,T}(r) \equiv g_a(r) \text { for } r \in\left[0, R_0\right] ; \\
g_{a,T}(r)>0\ \text{for}\ r\in\left(R_0,+\infty\right).
\end{array}\right.
$$
In order to investigate problem \eqref{eq1}, it is necessary to consider certain facts related to the existence of a normalized solution for the following problem:
\begin{equation}\label{eq2}
\left\{\begin{array}{l}
(-\Delta)^s u+\mu u=\lambda u+\left(I_\alpha *|u|^q\right)|u|^{q-2} u+\left(I_\alpha *|u|^p\right)|u|^{p-2} u, \quad x \in \mathbb{R}^N, \\
\int_{\mathbb{R}^N}|u|^2 \mathrm{d}x=a.
\end{array}\right.
\end{equation}
It is obvious that a solution $u$ to the problem \eqref{eq2} corresponds to a critical point of the  functional:
$$
I_\mu(u)=\frac{1}{2} \|u\|^2 +\frac{\mu}{2}\|u\|_2^2-\frac{1}{2q}Q(u)-\frac{1}{2p} P(u),
$$
restricted to  $S(a)$.
In this context, we can also consider the following functional:
$$
I_{\mu, T}(u):=\frac{1}{2} \|u\|^2+\frac{\mu}{2}\|u\|_2^2 -\frac{1}{2q}Q(u)-\frac{1}{2p} \tau(\|u\|)P(u).
$$
It is easy to see that $I_{\mu} \in C^1\left(H^s\left(\mathbb{R}^N\right), \mathbb{R}\right)$ and $I_{\mu, T} \in C^1\left(H^s\left(\mathbb{R}^N\right), \mathbb{R}\right)$ respectively. Moreover, according to \eqref{1.3}, the following inequality holds for any $u \in S(a)$:
\begin{equation}\label{2.4}
\begin{aligned}
I_{\mu, T}(u) &\geq \frac{1}{2}\|u\|^2-\frac{1}{2 q} C_{\alpha, q}\|u\|^{2q \gamma_{q, s}}a^{ q\left(1-\gamma_{q, s}\right)}-\frac{1}{2 p} C_{\alpha, p}\tau\left(\| u\|\right)\|u\|^{2 p \gamma_{p, s}}a^{ p\left(1-\gamma_{p, s}\right)} \\
&=g_{a,T}\left(\| u\|\right).
\end{aligned}
\end{equation}

\begin{lemma}\label{le2.1}
When $0<a_1 \leq a$, the functional $I_{\mu, T}$ exhibits a lower bound on the set $S\left(a_1\right)$.
\end{lemma}
{\bf Proof.} For any $u \in S\left(a_1\right)$, according to \eqref{2.4} and the properties of $g_{a,T}(r)$, it can be deduced that
\begin{equation*}
I_{\mu, T}(u)\geq g_{a_1,T}\left(\| u\|\right) \geq g_{a,T}\left(\| u\|\right) \geq \inf _{r \geq 0} g_{a,T}(r)>-\infty.
\end{equation*}
\hfill{$\Box$}

\begin{lemma} \label{le2.2}
Let $0<a_1 \leq a$, there exists $V_*>0$ such that if $\mu \in [0, V_*)$, there holds that \[\Upsilon_{\mu, T, a_1}:=\inf\limits _{u \in S\left(a_1\right)} I_{\mu, T}(u)<0.\]
\end{lemma}
{\bf Proof.}Let $u \in S\left(a_1\right)$ be fixed. For $t>0$, we define $u_t(x)=t^{\frac{N}{2}} u(t x)$. It follows that $u_t \in S\left(a_1\right)$ for all $t>0$. From $\tau \geq 0$, we can conclude that
\begin{equation}\label{2.5}
\begin{aligned}
I_{\mu, T}\left(u_t\right) & \leq \frac{1}{2} \|u_t\|^2 +\frac{\mu}{2}a_1-\frac{1}{2q} \int_{\mathbb{R}^N}\left(I_\alpha *|u_t|^q\right)|u_t|^q \\
& =\frac{1}{2} t^{2s} \|u\|^2+\frac{\mu}{2}a_1-\frac{1}{2q} t^{Nq-N-\alpha} Q(u).
\end{aligned}
\end{equation}
Due to $q <\frac{N+2s+\alpha}{N}$, we infer that $\frac{1}{2} t^{2s} \|u\|^2-\frac{1}{2q} t^{Nq-N-\alpha} Q(u)<0$ for sufficiently small $t>0$. Therefore, in accordance with \eqref{2.5}, there exists $V_*>0$ such that $I_{\mu, T}\left(u_t\right)<0$ for $\mu \in [0, V_*)$.
\hfill{$\Box$}

\begin{lemma} \label{le2.3}
For $\mu \in\left[0, V_*\right)$ and  $a_1 \in(0, a]$ and $u \in S\left(a_1\right)$ such that $I_{\mu, T}(u)<0$, then $\| u\|<R_0$. Furthermore, for all $v$ with $\|v\|^2_2 \leq a$ and close to $u$ in $H^s\left(\mathbb{R}^N\right)$, we have $I_{\mu, T}(v)=$ $I_\mu(v)$.
\end{lemma}
{\bf Proof.} From \eqref{2.4} and $I_{\mu, T}(u)<0$, we infer that
$
g_{a,T}\left(\| u\|\right)\leq g_{a_1,T}\left(\| u\|\right)\leq I_{\mu, T}(u)<0.
$
Consequently, it can be deduced from the properties of $g_{a,T}(r)$ that $\|u\|<R_0$.
By the continuity of $I_{\mu, T}$,  $I_{\mu, T}(v)<0$ for all $v$ close to $u$ in $H^s\left(\mathbb{R}^N\right)$. After further observing that   $\|v\|^2_2 \leq a$, it follows that $\| v\|<R_0$, and consequently we have $I_{\mu, T}(v)=I_\mu(v)$. \hfill{$\Box$}

For any $a_1 \in(0, a]$, we define
$$
m_\mu\left(a_1\right):=\inf _{u \in V\left(a_1\right)} I_\mu(u),
$$
where $V\left(a_1\right):=\left\{u \in S\left(a_1\right):\| u\|<R_0\right\}$.

\begin{lemma} \label{le2.4}
For  $\mu \in\left[0, V_*\right)$ and  $a_1 \in(0, a]$, it follows that
$
\Upsilon_{\mu, T, a_1}=m_\mu\left(a_1\right).
$
\end{lemma}
{\bf Proof.}
Based on the definition of $\Upsilon_{\mu, T, a_1}$ and Lemma \ref{le2.2}, we select a sequence $\{u_n\}\subset S(a_1)$ such that $I_{\mu, T}(u_n)\to\Upsilon_{\mu, T, a_1}<0$.  Consequently, for sufficiently large $n$, we achieve $I_{\mu, T}(u_n)<0$. Therefore, by Lemma \ref{le2.3}, $\|u_n\|<R_0$. This implies that $u_n$ belongs to $V(a_1)$. As a result, $m_\mu\left(a_1\right)\leq \Upsilon_{\mu, T, a_1}$. On the contrary, since $I_{\mu, T}(u) = I_\mu(u)$ in $V\left(a_1\right)$, we get that
$$
m_\mu\left(a_1\right)=\inf _{u \in V\left(a_1\right)} I_\mu(u)=\inf _{u \in V\left(a_1\right)} I_{\mu,T}(u)\geq\inf _{u \in S\left(a_1\right)} I_{\mu,T}(u)=\Upsilon_{\mu, T, a_1}.
$$
\hfill{$\Box$}

\begin{lemma}\label{le2.5}
Let $0<a_1\leq a$ and $r_1>0$ such that $w_{a_1}\left(r_1\right) \geq 0$, then for any $a_2 \in\left(0, a_1\right]$,
$$
w_{a_2}\left(r_2\right) \geq 0, \quad \text { for} \quad r_2 \in\left[\sqrt{\frac{a_2}{a_1}} r_1, r_1\right] .
$$
\end{lemma}
{\bf Proof.} Since $w_a(r)$ in \eqref{2.2} is a non-increasing function in terms of $a$, it is evident that
$
w_{a_2}\left(r_1\right) \geq w_{a_1}\left(r_1\right) \geq 0 .
$
Through direct calculations,
$$
w_{a_2}\left(\sqrt{\frac{a_2}{a_1}} r_1\right) \geq w_{a_1}\left(r_1\right) \geq 0 .
$$
Based on the above inequalities and the property of  $w_{a_2}(r)$, the conclusion is deemed valid.
\hfill{$\Box$}

\begin{lemma}\label{le2.6}
For $\mu \in\left[0, V_*\right)$, we have
\begin{itemize}
\item[$(1)$] The function $b \in(0, a] \mapsto m_\mu\left(b\right)$ exhibits continuity.
\item[$(2)$] Given $0<a_1< a_2 \leq a$, it follows that $\frac{a_1}{a_2} m_\mu\left(a_2\right)\leq m_\mu\left(a_1\right)$. Furthermore, if the value of $m_\mu$ at $a_1$ is attained, then the inequality becomes strict.
\end{itemize}
\end{lemma}
{\bf Proof.}
$(1)$ For any $b \in\left(0, a\right]$ and a sequence $\left\{a_n\right\} \subset\left(0, a\right)$ such that $ a_n \to b $ as $ n \to \infty$, we can refer to the definition of $ m_\mu(a_n)$. For each integer \( n\), there exists an element \( u_n \in V(a_n) \) satisfying:
\begin{equation}\label{2.6}
I_\mu\left(u_n\right)<m_\mu\left(a_n\right)+\frac{1}{n}
\end{equation}
and
\begin{equation}\label{2.7}
I_\mu\left(u_n\right)<0.
\end{equation}
Since it holds that $u_n \in V\left(a_n\right)$, we have $\|u_n\|_2^2=a_n$ and $\|u_n\|< R_0$. This indicates that the sequence  $\left\{u_n\right\}$  is bounded in $H^s\left(\mathbb{R}^N\right)$. Consequently, both $Q(u_n)$ and $P(u_n)$ are bounded. Following reasoning similar to what was done for \eqref{2.1}:
\begin{equation*}
 I_\mu(u_n)\geq g_{a_n}(\|u_n\|)=\frac{1}{2}\|u_n\|^2 w_{a_n}(\|u_n\|).
\end{equation*}
Combining this inequality with \eqref{2.7} leads us to conclude that $w_{a_n}(\|u_n\|)<0$. Thus, by applying Lemma \ref{le2.5} along with noting that   $w_a(R_0)=0$, we derive
\begin{equation}\label{2.8}
\|u_n\|<\sqrt{\frac{a_n}{a}}R_0.
\end{equation}
Next, define $v_n=\sqrt{\frac{b}{a_n}} u_{n}$ and then $v_{n}\in S(b)$. From the inequality derived above,
$$
\left\|v_n\right\|^2=\frac{b}{a_n}\left\|u_n\right\|^2<\frac{b}{a_n} \cdot \frac{a_n}{a} R_0^2\leq R_0^2.
$$
Thus, $v_n \in V(b)$. This, together with the boundedness of both terms involving   $Q(u_n)$ and $P(u_n)$ along with the convergence  of $a_n$ to $b$, 
\begin{equation*}
\begin{aligned}
m_\mu (b) &\leq I_\mu\left(v_n\right)=I_\mu\left(u_n\right)+\frac{1}{2}\left(\frac{b}{a_n}-1\right)\left\|u_n\right\|^2+\frac{\mu}{2}\left(\frac{b}{a_n}-1\right)\left\|u_n\right\|_2^2\\
&-\frac{1}{2 q}\left(\left(\frac{b}{a_n}\right)^q-1\right) Q\left(u_n\right)
 -\frac{1}{2p}\left(\left(\frac{b}{a_n}\right)^{p}-1\right) P\left(u_n\right) \\
= & I_\mu\left(u_n\right)+o_n(1) .
\end{aligned}
\end{equation*}
Combining this with \eqref{2.6}, we can conclude
\begin{equation}\label{2.9}
m_\mu(b) \leq \liminf_{n \rightarrow \infty} m_\mu\left(a_n\right) .
\end{equation}
On the other hand, assume there exists a minimizing sequence  $\left\{e_n\right\} \subset V(b)$  for $m_\mu (b)$ with $I_\mu\left(e_n\right)<0$. Let $f_n=\sqrt{\frac{a_n}{b}} e_n$. Thus,  $f_n\in S(a_n)$ follows. By employing methods akin to those used in proving \eqref{2.8}, one can show $f_n \in V\left(a_n\right)$ as well. Therefore,
$$
m_\mu\left(a_n\right) \leq I_\mu\left(f_n\right)=I_\mu\left(e_n\right)+o_n(1)=m_\mu(b)+o_n(1),
$$
leading us to conclude
\begin{equation}\label{2.10}
\limsup_{n \rightarrow \infty} m_\mu\left(a_n\right) \leq m_\mu(b).
\end{equation}
By combining results from \eqref{2.9} and \eqref{2.10}, it is evident that  $m_\mu\left(a_n\right) \rightarrow m_\mu(b)$ as $n \rightarrow \infty$.

$(2)$  Based on Lemma \ref{le2.4} and the definition of $m_\mu(a_1)$, for any sufficiently small $\varepsilon>0$, there exists an element $u \in V\left(a_1\right)$ such that
\begin{equation}\label{2.11}
I_\mu(u) \leq m_\mu\left(a_1\right)+\varepsilon \text { and } I_\mu(u)<0 .
\end{equation}
Following a similar argument as in the proof of  \eqref{2.8},
$
\|u\|^2<\frac{a_1}{a} R^2_0.
$
Define  $\delta=\frac{a_2}{a_1}>1$ and set $v(x)=u\left(\delta^{-\frac{1}{N}} x\right)$. Consequently, we have  $\|v\|_2^2=a_2$ and
$$
\|v\|^2=\delta^{1-\frac{2 s}{N}}\|u\|^2<\delta\|u\|^2<\frac{a_2}{a_1} \cdot \frac{a_1}{a} R^2_0<R^2_0.
$$
Thus, it follows that $v \in V\left(a_2\right)$. From this result along with \eqref{2.11}, we can conclude
$$
\begin{aligned}
m_\mu\left(a_2\right) & \leq I_\mu(v)=\frac{\delta^{1-\frac{2 s}{N}}}{2}\|u\|^2+\frac{1}{2}a_2-\frac{1}{2 q} \delta^{\frac{N+\alpha}{N}} T(u)-\frac{1}{2p} \delta^{\frac{N+\alpha}{N}} W(u) \\
& <\delta I_\mu(u) \leq \delta\left(m_\mu\left(a_1\right)+\varepsilon\right) .
\end{aligned}
$$
Taking the limit  $\varepsilon \rightarrow 0$ in this inequality gives us
$
m_\mu\left(a_2\right) \leq \frac{a_2}{a_1} m_\mu\left(a_1\right) .
$
Specifically, if $m_\mu\left(a_1\right)$ is attained, setting $\varepsilon=0$ results in
$
m_\mu\left(a_2\right)<\frac{a_2}{a_1} m_\mu\left(a_1\right) .
$
\hfill{$\Box$}

\begin{lemma} \label{le2.7}
For $\mu \in\left[0, V_*\right)$, we have
\begin{itemize}
\item[$(1)$]  $a_1 \in(0, a] \mapsto \Upsilon_{\mu, T, a_1}$ is continuous.
\item[$(2)$] Let $0<a_1<a_2 \leq a$, then $\frac{a_1}{a_2} \Upsilon_{\mu, T, a_2}\leq\Upsilon_{\mu, T, a_1}$. Moreover, if  $\Upsilon_{\mu, T, a_1}$ is reached, the inequality is strict.
\end{itemize}
\end{lemma}
{\bf Proof.} From Lemmas \ref{le2.4} and \ref{le2.6}, we can  get the reults.

\begin{lemma}\label{le2.8}
Let $\mu \in\left[0, V_*\right)$, $a_1 \in(0, a]$ and $\left\{u_n\right\} \subset S\left(a_1\right)$ be a minimizing sequence with respect to $\Upsilon_{\mu, T, a_1}$. Then, up to subsequence, either
\begin{itemize}
\item[$(I)$] $\left\{u_n\right\}$ strongly converges, or
\item[$(II)$] there exists $\left\{y_n\right\} \subset \mathbb{R}^N$ with $\left|y_n\right| \rightarrow \infty$ such that the sequence $v_n(x)=u_n\left(x+y_n\right)$ strongly converges to a function $v \in S\left(a_1\right)$ with $I_{\mu, T}(v)=\Upsilon_{\mu, T, a_1}$.
\end{itemize}
\end{lemma}
{\bf Proof.} According to Lemmas \ref{le2.2} and \ref{le2.3}, for sufficiently large $n$, $\left\| u_n\right\|<R_0$. It follows that there exists an element $u \in H^s\left(\mathbb{R}^N\right)$ such that, up to a subsequence, $u_n \rightharpoonup u$ in the space $H^s\left(\mathbb{R}^N\right)$ and $u_n \to u,\ \text{a.e. in} \ \mathbb{R}^N$. Now, let us consider the following three potential scenarios:

Case $1$. If $u \not \equiv 0$ and $\|u\|^2_2=b \neq a_1$, it follows that $b \in\left(0, a_1\right)$. Define $v_n=u_n-u$, then, by Br\'{e}zis-Lieb Lemma and Lemma 3.6 of \cite{meng-he}, we have
\begin{equation}\label{2.12}
\left\|u_n\right\|^2=\left\|v_n\right\|^2+\|u\|^2+o_n(1),
\end{equation}
\begin{equation}\label{2.13}
\left\|u_n\right\|_2^2=\left\|v_n\right\|_2^2+\|u\|_2^2+o_n(1),
\end{equation}
\begin{equation}\label{2.14}
Q(u_n)=Q(v_n)+Q(u)+o_n(1)\ \text{and}\ P(u_n)=P(v_n)+P(u)+o_n(1).
\end{equation}
Let $d_n=\left\|v_n\right\|^2_2$, without loss of generality, we assume $d_n\rightarrow d>0$ and by \eqref{2.13}, then $a_1=d+b$. Noting that $d_n \in\left(0, a_1\right)$ for $n$ large, \eqref{2.12}-\eqref{2.14}, and $\tau$ is non-increasing, we have
$$
\begin{aligned}
\Upsilon_{\mu, T, a_1}+o_n(1)=I_{\mu, T}\left(u_n\right)= & \frac{1}{2}\left\|v_n\right\|^2+\frac{\mu}{2}\left\|v_n\right\|_2^2-\frac{1}{2q}Q(v_n)-\frac{1}{2p} \tau\left(\left\| u_n\right\|\right)P(v_n) \\
& +\frac{1}{2}\| u\|^2+\frac{\mu}{2}\left\|u\right\|_2^2-\frac{1}{2q}Q(u)-\frac{1}{2p} \tau\left(\left\|u_n\right\|\right)P(u)+o_n(1) \\
\geq & I_{\mu, T}\left(v_n\right)+I_{\mu, T}(u)+o_n(1) \\
\geq & \Upsilon_{\mu, T, d_n}+I_{\mu, T}(u)+o_n(1).
\end{aligned}
$$
Letting $n \rightarrow+\infty$, we find that
\begin{equation}\label{2.16}
\Upsilon_{\mu, T, a_1}\geq  \Upsilon_{\mu, T, d}+I_{\mu, T}(u).
\end{equation}
If $I_{\mu, T}(u)>\Upsilon_{\mu, T, b}$, by \eqref{2.16} and Lemma \ref{le2.7}, we can obtain
$$
\begin{aligned}
\Upsilon_{\mu, T, a_1} & > \Upsilon_{\mu, T, d}+\Upsilon_{\mu, T, b}\geq \frac{d}{a_1} \Upsilon_{\mu, T, a_1}+\frac{b}{a_1}\Upsilon_{\mu, T, a_1}=\Upsilon_{\mu, T, a_1},
\end{aligned}
$$
which is a contradiction. Thus, $\Upsilon_{\mu, T, b}$ can be reached. Then using \eqref{2.16} and  Lemma \ref{le2.7}, we have
$$
\begin{aligned}
\Upsilon_{\mu, T, a_1} & \geq \frac{d}{a_1} \Upsilon_{\mu, T, a_1}+\Upsilon_{\mu, T, b}>\frac{d}{a_1} \Upsilon_{\mu, T, a_1}+\frac{b}{a_1} \Upsilon_{\mu, T, a_1}=\Upsilon_{\mu, T, a_1},
\end{aligned}
$$
which is a contradiction. So this possibility cannot exist.

Case $2$. If $\|u\|_2=a_1$, then $u_n \rightarrow u$ in $L^2\left(\mathbb{R}^N\right)$ and consequently,
\begin{equation}\label{n2.16}
u_n \rightarrow u \ \text{in}\  L^t\left(\mathbb{R}^N\right),\quad  \text{for all}\ t \in\left(2,2_s^*\right).
\end{equation}

For $p<2_{\alpha,s}^*$, in accordance with \eqref{1.3} and \eqref{n2.16}, it follows that $Q(u_n)\to Q(u)$ and $P(u_n)\to P(u)$. Thanks to the lower semi-continuity of the norm,
$$
\begin{aligned}
\Upsilon_{\mu, T, a_1}=\lim_{n\to\infty}I_{\mu, T}\left(u_n\right)
&= \lim_{n\to\infty}I_\mu\left(u_n\right)\\
&=\lim_{n\to\infty}\left[\frac{1}{2} \|u_n\|^2+\frac{\mu}{2} \|u_n\|_2^2-\frac{1}{2q} Q(u_n)-\frac{1}{2p}P(u_n)\right] \\
& \geq I_\mu(u)=I_{\mu,T}(u)\geq\Upsilon_{\mu, T, a_1}.
\end{aligned}
$$
Therefore, $\|u_n\|\to\|u\|$. Consequently, $u_n\to u$ in $H^s\left(\mathbb{R}^N\right)$. This implies that $(I)$ occurs.

For $p=2_{\alpha,s}^*$, given $\left\| v_n\right\| \leq\left\| u_n\right\|<R_0$ for sufficiently large $n$, $T(v_n)=o_n(1)$ and  \eqref{1.3},
\begin{equation}\label{2.17}
\begin{aligned}
I_{\mu, T}\left(v_n\right)= I_\mu\left(v_n\right) &\geq\frac{1}{2} \|v_n\|^2-\frac{1}{2q} Q(v_n)-\frac{1}{2p}P(v_n) \\
& \geq \frac{1}{2} \|v_n\|^2-\frac{1}{2p}P(v_n)+o_n(1) \\
& \geq\frac{1}{2} \|v_n\|^2\left(1-\frac{1}{2_{\alpha,s}^*} \frac{1}{S_\alpha^{2_{\alpha,s}^*}}\left\| v_n\right\|^{22_{\alpha,s}^*-2}\right)+o_n(1) \\
& \geq\frac{1}{2} \|v_n\|^2\left(1-\frac{1}{2_{\alpha,s}^*} \frac{1}{S_\alpha^{2_{\alpha,s}^*}}R_0^{22_{\alpha,s}^*-2}\right)+o_n(1) \\
& =\frac{1}{2} \|v_n\|^2 \frac{1}{q} C_{\alpha, q} a^{q\left(1-\gamma_{q,s}\right)} R_0^{2q \gamma_q-2}+o_n(1),
\end{aligned}
\end{equation}
where we have used  $w_a\left(R_0\right)=1-\frac{1}{2_{\alpha,s}^*} \frac{1}{S_\alpha^{2_{\alpha,s}^*}}R_0^{22_{\alpha,s}^*-2}-\frac{1}{q} C_{\alpha, q} a^{q\left(1-\gamma_{q,s}\right)} R_0^{2q \gamma_q-2}=0$ in the last equality in \eqref{2.17}. Now, it is important to recall that
\begin{equation}\label{2.18}
 I_{\mu, T}\left(u_n\right) \geq I_{\mu, T}\left(v_n\right)+I_{\mu, T}(u)+o_n(1).
\end{equation}
For any $u \in S\left(a_1\right)$, we have $I_{\mu, T}(u) \geq \Upsilon_{\mu, T,a_1}$. Combining this with \eqref{2.17} and \eqref{2.18}, 
\begin{equation*}
 I_{\mu, T}\left(u_n\right) \geq \frac{1}{2} \|v_n\|^2 \frac{1}{q} C_{\alpha, q} a^{q\left(1-\gamma_{q,s}\right)} R_0^{2q \gamma_q-2}+\Upsilon_{\mu, T,a_1}+o_n(1).
\end{equation*}
Noting $I_{\mu, T}\left(u_n\right)\to \Upsilon_{\mu, T,a_1}$ and letting $n\to\infty$ in the above inequality, we can get
$\left\| v_n\right\|^2 \rightarrow 0$, and consequently $u_n \rightarrow u$ in $H^s\left(\mathbb{R}^N\right)$. This implies that $(I)$ occurs.

Case $3$.  If $u \equiv 0$, i.e., $u_n \rightharpoonup 0$ in $H^s\left(\mathbb{R}^N\right)$, then we assert the existence of $R, \beta>0$, and $\left\{y_n\right\} \subset \mathbb{R}^N$ such that
\begin{equation}\label{nne2.18}
\int_{B_R\left(y_n\right)}\left|u_n\right|^2 \mathrm{~d} x \geq \beta, \quad \text { for all } n .
\end{equation}
Otherwise, we would have $u_n \rightarrow 0$ in $L^t\left(\mathbb{R}^N\right)$ for all $t \in\left(2,2_s^*\right)$. Consequently,  for  $p<2_{\alpha,s}^*$, according to \eqref{1.3},  $Q(u_n)=o_n(1)$ and $P(u_n)=o_n(1)$. Thus, we have
$$
I_{\mu, T}\left(u_n\right) = \frac{1}{2}\left\| u_n\right\|^2+\frac{\mu}{2}\left\| u_n\right\|_2^2+o_n(1),
$$
which contradicts $I_{\mu, T}\left(u_n\right) \rightarrow \Upsilon_{\mu, T, a_1}<0$. For  $p=2_{\alpha,s}^*$, similarly to \eqref{2.17}, we obtain that
$$
I_{\mu, T}\left(u_n\right) \geq\frac{1}{2} \|u_n\|^2 \frac{1}{q} C_{\alpha, q} a^{q\left(1-\gamma_{q,s}\right)} R_0^{2q \gamma_q-2}+o_n(1).
$$
This leads to a contradiction. Therefore, \eqref{nne2.18} holds and it is obvious that $\left|y_n\right| \rightarrow +\infty$. By defining $v_n(x)=u_n\left(x+y_n\right)$, it is evident that the sequence $\left\{v_n\right\} \subset S\left(a_1\right)$ also serves as a minimizing sequence for $\Upsilon_{\mu, T, a_1}$. Additionally, there exists a non-zero function $v \in H^s(\mathbb{R}^N) $ such that $v_n\rightharpoonup v $ in   $ H^s(\mathbb{R}^N)$. Using similar arguments as those applied in the first two scenarios of this proof, we can deduce that $v_n \rightarrow v$ in $H^s(\mathbb{R}^N)$, which indicates that $(II)$ takes place.
\hfill{$\Box$}

\begin{lemma} \label{le2.9}
Let  $\mu \in\left(0, V_*\right]$ and $a_1 \in(0, a]$, the value of $\Upsilon_{\mu, T, a_1}$ can be achieved.
\end{lemma}
{\bf Proof.} According to the definition of $\Upsilon_{\mu, T, a_1}$, there exists a  minimizing sequence $\left\{u_n\right\} \subset S\left(a_1\right)$ such that $I_{\mu, T}\left(u_n\right) \rightarrow \Upsilon_{\mu, T, a_1}$ as $n \rightarrow +\infty$. By applying Lemma \ref{le2.8}, it follows that there exists $u \in S\left(a_1\right)$ satisfying $I_{\mu, T}(u)=\Upsilon_{\mu, T, a_1}$.
\hfill{$\Box$}

\begin{lemma}\label{le2.10}
Let $a_1 \in(0, a]$ be fixed and suppose $0<\mu_1<\mu_2 \leq V_*$, we have $\Upsilon_{\mu_1, T, a_1}<\Upsilon_{\mu_2, T, a_1}$.
\end{lemma}
{\bf Proof.} Applying Lemma \ref{le2.9}, there exists $u \in S\left(a_1\right)$ such that $I_{\mu_2, T}(u)=\Upsilon_{\mu_2, T, a_1}$. Then, it follows that $\Upsilon_{\mu_1, T, a_1}\leq I_{\mu_1, T}(u)<  I_{\mu_2, T}(u)=\Upsilon_{\mu_2, T, a_1}$.
\hfill{$\Box$}

\section{The nonautonomous case}

From this point onward, we assume that $\|V\|_\infty<V_*$, with $V_*$ being the value derived in Lemma  \ref{le2.2}.  The notations $I_{0, T}$, $I_{\infty, T}$ : $H^s\left(\mathbb{R}^N\right) \rightarrow \mathbb{R}$ will be used to represent the following functionals:
$$
I_{0, T}(u):=\frac{1}{2}\| u\|^2-\frac{1}{2q}Q(u)-\frac{\tau\left(\| u\|\right)}{2p}P(u)
$$
and
$$
I_{\infty, T}(u):=\frac{1}{2}\|u\|^2+\frac{1}{2} \int_{\mathbb{R}^N} V_{\infty}|u|^2 d x-\frac{1}{2q}Q(u)-\frac{\tau\left(\| u\|\right)}{2p}P(u).
$$
Define
$$
\Upsilon_{0, T, a}:=\inf _{u \in S(a)} I_{0, T}(u) \quad \text{and}\quad \Upsilon_{\infty, T, a}:=\inf _{u \in S(a)} I_{\infty, T}(u),
$$
respectively. By Lemma \ref{le2.9},   there exist $u_0, u_{\infty} \in S(a)$ such that $I_{0, T}\left(u_0\right)=\Upsilon_{0, T, a}$ and $I_{\infty, T}\left(u_{\infty}\right)=\Upsilon_{\infty, T, a}$. Furthermore, it follows from Lemmas \ref{le2.3} and  \ref{le2.10} that $\Upsilon_{0, T, a}<\Upsilon_{\infty, T, a}<0$. Similarly to the proof of Lemma \ref{le2.3}, the following result can be obtained.

\begin{lemma}\label{le3.2}
 Let $a_1 \in(0, a]$ and $u \in S\left(a_1\right)$ such that $J_{\epsilon, T}(u) < 0$, then $\| u\|<R_0$. Furthermore, for all $v$ satisfying $\|v\|^2_2 \leq a$ and close to $u$ in $H^s\left(\mathbb{R}^N\right)$, $J_{\epsilon, T}(v)=J_{\epsilon}(v)$.
\end{lemma}

\begin{lemma}\label{le3.3}
Let $\left\{u_n\right\} \subset S(a)$ with $J_{\epsilon, T}\left(u_n\right)\to m<\Upsilon_{\infty, T, a}$. If $u_n \rightharpoonup u_{\epsilon}$ in $ H^s\left(\mathbb{R}^N\right)$, then $u_{\epsilon}\not\equiv0$.
\end{lemma}
{\bf Proof.}Assume for the sake of contradiction that $u_{\epsilon}\equiv 0$, then
\[
m+o_n(1)=J_{\epsilon, T}\left(u_n\right)=I_{\infty, T}\left(u_n\right)+\frac{1}{2} \int_{\mathbb{R}^N}\left(V(\epsilon x)-V_{\infty}\right)\left|u_n\right|^2 \mathrm{d} x .
\]
According to condition $(V_2)$, for any given $\xi > 0 $, there exists  $R > 0$ such that
\[
V(x) \geq V_{\infty} - \xi, \quad |x| \geq R.
\]
Thus, we have
\[
m + o_n(1) = J_{\epsilon, T}(u_n) \geq I_{\infty, T}(u_n) + \frac{1}{2} \int_{B_{R/\epsilon}(0)} (V(\epsilon x)-V_\infty)|u_n|^2 \mathrm{d}x -\frac{\xi}{2}\int_{B_{R/\epsilon}^{c}(0)} |u_n|^2 \mathrm{d}x.
\]
By applying Lemma \ref{le3.2}, $\left\{u_n\right\}$ is bounded in $ H^s(\mathbb{R}^N)$. It follows from $u_n \rightarrow 0$ in $L^2(B_{R/\epsilon}(0))$ that
$
m + o_n(1)\geq I_{\infty,T}(u_n)-C\xi,
$
which implies
$
m+o_n(1)\geq \Upsilon _{\infty,T,a}-C\xi.
$
Since the choice of $\xi>0$ is arbitrary, $m \geq \Upsilon _{\infty,T,a}$, which is a contradiction. Then, $u_\epsilon\not\equiv0$.
\hfill{$\Box$}

\begin{lemma}\label{le3.4}
Assume that \eqref{1.5} holds true if $p=2_{\alpha,s}^*$ and let $\left\{u_n\right\}$ be a $(P S)_c$ sequence of $J_{\epsilon, T}$ restricted to $S(a)$ with $c<\Upsilon_{\infty, T, a}$ and $u_n \rightharpoonup u_\epsilon$ in $H^s\left(\mathbb{R}^N\right)$. If $u_n$ does not converge to $u_\epsilon$ strongly in $H^s\left(\mathbb{R}^N\right)$, then there exist constants $\epsilon_0>0$ and $\beta>0$, independent of $\epsilon$, such that for $\epsilon \in\left(0, \epsilon_0\right)$,
$$
\limsup _{n \rightarrow+\infty}\left\|u_n-u_\epsilon\right\|_2 \geq \beta.
$$
\end{lemma}
{\bf Proof.} By Lemma \ref{le3.2}, $\left\| u_n\right\|<R_0$ for sufficiently large $n$. Consequently,
$$
J_\epsilon\left(u_n\right) \rightarrow c \quad \text { and } \quad J_\epsilon|_{S(a)} ^{\prime}\left(u_n\right) \rightarrow 0 \quad \text { as } n \rightarrow+\infty .
$$
Introducing the functional $\pi: H^s(\mathbb{R}^N) \rightarrow \mathbb{R}$ defined by:
$$
\pi(u)=\frac{1}{2} \int_{\mathbb{R}^N}|u|^2 \mathrm{d}x,
$$
we have $S(a)=\pi^{-1}\left(a / 2\right)$. According to \cite[Proposition 5.12]{MR1400007}, there exists $\left\{\lambda_n\right\} \subset \mathbb{R}$ such that
\begin{equation}\label{3.1}
\left\|J_\epsilon^{\prime}\left(u_n\right)-\lambda_n \pi^{\prime}\left(u_n\right)\right\|_{H^{-s}\left(\mathbb{R}^N\right)} \rightarrow 0, \quad \text { as } n \rightarrow+\infty .
\end{equation}
Due to the boundedness of $\left\{u_n\right\}$ in $H^s\left(\mathbb{R}^N\right)$, $\left\{\lambda_n\right\}$ is also bounded. Then, up to subsequence, $\lambda_n\rightarrow\lambda_\epsilon$ as $n \rightarrow +\infty$. This, combined with \eqref{3.1}, leads to
\begin{equation}\label{3.2}
J_\epsilon^{\prime}\left(u_\epsilon\right)-\lambda_\epsilon \pi^{\prime}\left(u_\epsilon\right)=0, \quad \text { in } \ H^{-s}\left(\mathbb{R}^N\right).
\end{equation}
Set $v_n:=u_n-u_\epsilon$, then it follows from \eqref{3.1} and \eqref{3.2} that
\begin{equation}\label{3.3}
\left\|J_\epsilon^{\prime}\left(v_n\right)-\lambda_\epsilon \pi^{\prime}\left(v_n\right)\right\|_{H^{-s}\left(\mathbb{R}^N\right)} \rightarrow 0, \quad \text { as } n \rightarrow+\infty,
\end{equation}
By $J_\epsilon\left(u_n\right) \rightarrow c<\Upsilon_{\infty, T, a}$, $\lambda_n\to\lambda_\epsilon$ and \eqref{3.1}, it is determined that
$$
\begin{aligned}
& \Upsilon_{\infty, T, a}> \lim _{n \rightarrow+\infty} J_\epsilon\left(u_n\right) \\
&=\lim _{n \rightarrow+\infty}\left(J_\epsilon\left(u_n\right)-\frac{1}{2} J_\epsilon^{\prime}\left(u_n\right) u_n+\frac{1}{2} \lambda_n\left\|u_n\right\|_2^2+o_n(1)\right) \\
&=\lim _{n \rightarrow+\infty}\left[\left(\frac{1}{2}-\frac{1}{2q}\right) Q(u_n)+\left(\frac{1}{2}-\frac{1}{2p}\right)  P(u_n)+\frac{1}{2} \lambda_n a+o_n(1)\right] \\
& \geq \frac{1}{2} \lambda_\epsilon a,
\end{aligned}
$$
which implies that
\begin{equation}\label{3.4}
\lambda_\epsilon \leq \frac{2 \Upsilon_{\infty, T, a}}{a}<0.
\end{equation}
By referring to \eqref{3.3}, we have
\begin{equation}\label{3.5}
\left\|v_n\right\|^2+\int_{\mathbb{R}^N} V(\epsilon x)\left|v_n\right|^2 \mathrm{d} x-\lambda_\epsilon\left\|v_n\right\|_2^2=Q(v_n)+P(v_n)+o_n(1),
\end{equation}
which, combined with \eqref{3.4},  yields
\begin{equation}\label{3.6}
\left\| v_n\right\|^2+\int_{\mathbb{R}^N} V(\epsilon x)\left|v_n\right|^2 \mathrm{d} x-\frac{2 \Upsilon_{\infty, T, a}}{a}\left\|v_n\right\|_2^2 \leq Q(v_n)+P(v_n)+o_n(1) .
\end{equation}
By using \eqref{3.6} and \eqref{1.3}, there exists $C_1$ independent of $\epsilon$ such that
\begin{equation}\label{3.7}
C_1\left\|v_n\right\|_{H^{s}\left(\mathbb{R}^N\right)}^2\leq  C_{\alpha,q}\left\|v_n\right\|_{H^{s}\left(\mathbb{R}^N\right)}^{2q}+C_{\alpha,p}\left\|v_n\right\|_{H^{s}\left(\mathbb{R}^N\right)}^{2p}+o_n(1) .
\end{equation}
Since $v_n \nrightarrow 0$ in $H^s\left(\mathbb{R}^N\right)$, by the second inequality in \eqref{3.7}, there exists $C_2$ independent of $\epsilon$ such that $\left\|v_n\right\|_{H^{s}\left(\mathbb{R}^N\right)} \geq C_2$. Furthermore, according to the first inequality in \eqref{3.7}, 
\begin{equation}\label{3.8}
\limsup _{n \rightarrow+\infty}\left(Q(v_n)+P(v_n)\right) \geq C_3
\end{equation}
for some $C_3$ independent of $\epsilon$.

{\bf Case 1:} $p<2_{\alpha,s}^*$. Employing \eqref{1.3},
we arrive at
\begin{equation}\label{3.9}
\limsup _{n \rightarrow+\infty}Q(v_n) \leq C_{\alpha,q}\left(\limsup _{n \rightarrow+\infty}\|v_n\|_2\right)^{2q\left(1-\gamma_{q,s}\right)} R_0^{2q\gamma_{q,s}},
\end{equation}
and
\begin{equation}\label{3.10}
\limsup _{n \rightarrow+\infty}P(v_n) \leq C_{\alpha,p}\left(\limsup _{n \rightarrow+\infty}\|v_n\|_2\right)^{2p\left(1-\gamma_{p,s}\right)} R_0^{2p\gamma_{p,s}}.
\end{equation}
 Now the lemma follows from  \eqref{3.8}-\eqref{3.10}.

{\bf Case 2:} $p=2_{\alpha,s}^*$. We assert that there exists $\epsilon_0>0$ such that for $\epsilon\in(0,\epsilon_0)$,
$
\limsup _{n \rightarrow+\infty} Q(v_n) \geq C
$
for some $C>0$ independent of $\epsilon$. If this assertion holds, according to  \eqref{3.9}, we can obtain the result.
Otherwise, there exists a sequence $\left\{\epsilon_m\right\}$, such that
\begin{equation}\label{3.11}
\lim_{\epsilon_m \rightarrow 0^{+}} \limsup _{n \rightarrow+\infty}Q(v_n)=0.
\end{equation}
This, together with \eqref{3.8}, leads to
\begin{equation}\label{3.12}
\lim_{\epsilon_m\rightarrow 0^{+}} \limsup _{n \rightarrow+\infty} P(v_n) \geq C.
\end{equation}
It follows from \eqref{3.5} and \eqref{3.11} that
$$
\lim_{\epsilon_m \rightarrow 0^{+}} \limsup _{n \rightarrow+\infty} \|v_n\|^2\leq\lim_{\epsilon_m \rightarrow 0^{+}} \limsup _{n \rightarrow+\infty}  P(v_n).
$$
Applying \eqref{1.3} to the above inequality, we obtain that
$$
\lim_{\epsilon_m \rightarrow 0^{+}}\limsup _{n \rightarrow+\infty}\left\|v_n\right\|^2\leq \lim_{\epsilon_m \rightarrow 0^{+}}\limsup _{n \rightarrow+\infty}   \frac{1}{\mathcal{S}_\alpha^{2_{\alpha,s}^*}}\left\| v_n\right\|^{22_{\alpha,s}^*},
$$
which, combined with \eqref{3.12} and  $\left\| v_n\right\| \leq\left\| u_n\right\|<R_0$ for sufficiently large $n$, implies that
\begin{equation}\label{3.13}
R_0 \geq \lim_{\epsilon_m \rightarrow 0^{+}} \limsup _{n \rightarrow+\infty}\left\|v_n\right\| \geq  \mathcal{S}_\alpha^{\frac{N+\alpha}{2(\alpha+2s)}}.
\end{equation}
On the other hand, by \eqref{1.5} and the expression of $r_0$, we can get
$
r_0< \mathcal{S}_\alpha^{\frac{N+\alpha}{2(\alpha+2s)}}.
$
This, together with $R_0<r_0$, implies that
$
R_0< \mathcal{S}_\alpha^{\frac{N+\alpha}{2(\alpha+2s)}},
$
which  contradicts \eqref{3.13}.  The proof is complete.
\hfill{$\Box$}

\begin{lemma} \label{le3.5}
For $\epsilon \in (0, \epsilon_0)$ and let $\beta$ be given in Lemma \ref{le3.4}. Assume \eqref{1.5} holds true if $p=2_{\alpha,s}^*$ and 
$$
0 < \rho_0 \leq \frac{\beta}{a}(\Upsilon_{\infty, T, a}-\Upsilon_{0, T, a}).
$$
If $c<\Upsilon_{0,T,a}+\rho_0$, then the functional $J_{\epsilon,T}$ meets the $(PS)_c$ condition when restricted to $S(a)$.
\end{lemma}
{\bf Proof.}
Let $\{u_n\}$ be a $(PS)_c$ sequence of $J_{\epsilon, T}$ restricted to $S(a)$. Thanks to $c < \Upsilon_{\infty, T, a} < 0$, by Lemma \ref{3.2} we know $\{u_n\}$ is bounded in $H^s(\mathbb{R}^N)$. Consequently, there exists $u_\epsilon \in H^s(\mathbb{R}^N)$ such that $u_n \rightharpoonup u_\epsilon$. By Lemma \ref{le3.3}, $u_\epsilon\not\equiv0$. Define $v_n = u_n - u_\epsilon $. If $ u_n \to u_\epsilon $ strongly in $ H^s(\mathbb{R}^N) $, then we are done. However, if $ u_n $ does not converge strongly to $u_\epsilon $ in $ H^s(\mathbb{R}^N) $ for some $\epsilon \in\left(0, \epsilon_0\right)$, we refer to Lemma \ref{3.4}, which gives us
$
\limsup _{n \rightarrow+\infty}\left\|v_n\right\|^2_2 \geq \beta.
$
Letting $b=\left\|u_\epsilon\right\|^2_2$ and defining $d_n=\left\|v_n\right\|^2_2$, assume that as $n\to\infty$, we have $d_n\rightarrow d$.
This leads us to conclude that $d \geq \beta>0$ and $a=b+d$. Since for sufficiently large $n$ it holds true that $d_{n} \in(0,a)$,  we derive that
\begin{equation}\label{3.14}
c+o_n(1)=J_{\epsilon, T}\left(u_n\right) \geq J_{\epsilon, T}\left(v_n\right)+J_{\epsilon, T}\left(u_\epsilon\right)+o_n(1) .
\end{equation}
Given $v_n \rightharpoonup 0$ in $H^s\left(\mathbb{R}^N\right)$, similar reasoning as used in proving Lemma \ref{3.3} allows us to assert
\begin{equation}\label{3.15}
J_{\epsilon, T}\left(v_n\right) \geq I_{\infty, T}\left(v_n\right)-\delta C+o_n(1)
\end{equation}
for any $\delta>0$, where $C>0$ is a constant independent of $\delta, \epsilon$, and $n$. Thanks to \eqref{3.14} and \eqref{3.15},
\begin{equation*}
\begin{aligned}
c+o_n(1)=J_{\epsilon, T}\left(u_n\right) & \geq I_{\infty, T}\left(v_n\right)+J_{\epsilon, T}\left(u_\epsilon\right)-\delta C+o_n(1) \\
& \geq \Upsilon_{\infty, T, d_n}+\Upsilon_{0 , T, b}-\delta C+o_n(1).
\end{aligned}
\end{equation*}
Taking the limit as $n \rightarrow+\infty$ and utilizing Lemma \ref{le2.7} along with the arbitrary choice of $\delta > 0$,
$$
\begin{aligned}
c \geq \Upsilon_{\infty, T, d}+\Upsilon_{0, T, b} & \geq \frac{d}{a} \Upsilon_{\infty, T, a}+\frac{b}{a} \Upsilon_{0, T, a} \\
& =\Upsilon_{0, T, a}+\frac{d}{a}\left(\Upsilon_{\infty, T, a}-\Upsilon_{0, T, a}\right) \\
& \geq \Upsilon_{0, T, a}+\frac{\beta}{a}\left(\Upsilon_{\infty, T, a}-\Upsilon_{0, T, a}\right).
\end{aligned}
$$
This contradicts $c<\Upsilon_{0, T, a}+\frac{\beta}{a}\left(\Upsilon_{\infty, T, a}-\Upsilon_{0, T, a}\right)$. Therefore, $u_n \rightarrow u_\epsilon$ in $H^s\left(\mathbb{R}^N\right)$.
\hfill{$\Box$}

\section{Multiplicity result.}

Let $\delta>0$ be fixed and let $w$ satisfy
$I_{0,T}(w)=\Upsilon_{0, T, a}$. Consider a smooth non-increasing cut-off function $\eta$ defined on the interval $[0, \infty)$ such that $\eta(s)=1$ if $0 \leq s \leq \frac{\delta}{2}$ and $\eta(s)=0$ if $s \geq \delta$. For any point $y \in \mathcal{M}$, we can define
$$
\Psi_{\epsilon, y}(x)=\eta(|\epsilon x-y|) w((\epsilon x-y) / \epsilon),
$$
and
$$
\tilde{\Psi}_{\epsilon, y}(x)=\sqrt{a} \frac{\Psi_{\epsilon, y}(x)}{\|\Psi_{\epsilon, y}\|_2}.
$$
Through direct computation, we have
$\Psi_{\epsilon, y}(x)\to w(x)$
in $H^s\left(\mathbb{R}^N\right)$ and $\tilde{\Psi}_{\epsilon, y}(x)\to w(x)$ in $H^s\left(\mathbb{R}^N\right)$.
Define the mapping $\Phi_\epsilon: \mathcal{M} \rightarrow S(a)$ by setting  $\Phi_\epsilon(y)=\tilde{\Psi}_{\epsilon, y}$.

\begin{lemma}\label{le4.1}
$$
\lim _{\epsilon \rightarrow 0} J_{\epsilon,T}\left(\Phi_\epsilon(y)\right)=\Upsilon_{0, T, a}, \text { uniformly in } \quad y \in \mathcal{M}.
$$
\end{lemma}
{\bf Proof.} Noting  $\tilde{\Psi}_{\epsilon, y}(x)\to w(x)$ in $H^s\left(\mathbb{R}^N\right)$ and $V(0)=0$, we can easily derive this result.
\hfill{$\Box$}

For any $\delta>0$, let $R=R(\delta)>0$ be chosen such that $\mathcal{M}_\delta \subset B_R(0)$. Consider the mapping $\chi: \mathbb{R}^N \rightarrow \mathbb{R}^N$ defined as follows: for $|x| \leq R$, $\chi(x)=x$, and for $|x| \geq R$, $\chi(x)=\frac{R x}{|x|}$. Furthermore, let us define the function $\beta_\epsilon: S(a) \rightarrow \mathbb{R}^N$ by
$$
\beta_\epsilon(u)=\frac{\int_{\mathbb{R}^N} \chi(\epsilon x)|u|^2 \mathrm{d} x}{a}.
$$
According to \cite[Lemma 4.2]{alves-thin}, the following lemma can be derived, with the proof details omitted.

\begin{lemma}\label{le4.2}
$$
\lim _{\epsilon \rightarrow 0} \beta_\epsilon\left(\Phi_\epsilon(y)\right)=y \text {, uniformly in } \quad y \in \mathcal{M}.
$$
\end{lemma}

\begin{lemma} \label{le4.3}
Let $\epsilon_n \rightarrow 0$ and $\{u_n\} \subset S(a)$ such that $J_{\epsilon_n, T}(u_n)\to\Upsilon_{0, T, a}$. Then, there exists $\{\tilde{y}_n\} \subset \mathbb{R}^N$ such that $v_n(x) = u_n(x + \tilde{y}_n)$ has a subsequence that converges in $H^s(\mathbb{R}^N)$. Furthermore, after selecting an appropriate subsequence, $y_n:= \epsilon_n \tilde{y}_n$ converging to some limit $y$ within the set $\mathcal{M}$.
\end{lemma}
{\bf Proof.} It follows from 
$J_{\epsilon_n, T}\left(u_n\right)\geq I_{0, T}\left(u_n\right)\geq\Upsilon_{0, T, a},$ that
$I_{0, T}\left(u_n\right)\to\Upsilon_{0, T, a}.$
By applying Lemma  \ref{le2.8}, we find that $v_n \rightarrow v$ in $H^s\left(\mathbb{R}^N\right)$ and $v \in S(a)$.
Now, we will demonstrate that $\left\{y_n\right\}$ is bounded.
Suppose there exists a subsequence such that $\left|y_n\right| \rightarrow +\infty$, then
$$
\Upsilon_{0, T, a}=\lim _{n \rightarrow+\infty}\left(\frac{1}{2}\|v_n\|^2+\frac{1}{2}\int_{\mathbb{R}^N} V\left(\epsilon_n x+y_n\right)\left|v_n\right|^2 \mathrm{d} x-Q\left(v_n\right)-\tau(\|v_n\|)P\left(v_n\right) \right)
$$
which yields that
$$
\Upsilon_{0, T, a}=\frac{1}{2}\|v\|^2+\frac{1}{2}\int_{\mathbb{R}^N} V_
\infty\left|v\right|^2 \mathrm{d} x-Q\left(v\right)-\tau(v)P\left(v\right)\geq\Upsilon_{\infty, T, a}.
$$
This contradicts Lemma  \ref{le2.10}. Then, we assume that $y_n \rightarrow y$ in $\mathbb{R}^N$ and get similarly as above
$$
\Upsilon_{0, T, a}=\frac{1}{2}\|v\|^2+\frac{1}{2}\int_{\mathbb{R}^N} V(y)\left|v\right|^2 \mathrm{d} x-Q\left(v\right)-\tau(v)P\left(v\right) \geq\Upsilon_{V(y), T, a},
$$
which together with  Lemma \ref{le2.10}  implies that $V(y) \leq 0$. It follows from $(V_2)$ that $V(y)=0$ and $y \in \mathcal{M}$.\hfill{$\Box$}

Let $h:[0,+\infty) \rightarrow[0,+\infty)$ be a  function such that $h(\epsilon) \rightarrow 0$ as $\epsilon \rightarrow 0$. Consider the set
$$
\tilde{S}(a)=\left\{u \in S(a): J_{\epsilon, T}(u) \leq \Upsilon_{0, T, a}+h(\epsilon)\right\} .
$$
By Lemma \ref{le4.1}, the function $h(\epsilon)=\sup\limits _{y \in \mathcal{M}}\left|J_{\epsilon, T}\left(\Phi_\epsilon(y)\right)-\Upsilon_{0, T, a}\right|$ satisfies $h(\epsilon) \rightarrow 0$ as $\epsilon \rightarrow 0$.  Consequently, we have  $\Phi_\epsilon(y) \in \tilde{S}(a)$ for all $y \in \mathcal{M}$.

\begin{lemma}\label{le4.4}
For any $\delta>0$,
$$
\lim _{\epsilon \rightarrow 0} \sup _{u \in \tilde{S}(a)} \inf _{z \in \mathcal{M}_\delta}\left|\beta_\epsilon(u)-z\right|=0 .
$$
\end{lemma}
{\bf Proof.} Let $\epsilon_{n} \rightarrow 0$ and consider $u_{n} \in \tilde{S}(a)$ such that
$$
\inf _{z \in \mathcal{M}_\delta}\left|\beta_{\epsilon_n}\left(u_n\right)-z\right|=\sup _{u \in \tilde{S}(a)} \inf _{z \in \mathcal{M}_\delta}\left|\beta_{\epsilon_n}\left(u_n\right)-z\right|+o_n(1) .
$$
It suffices to identify a sequence $\left\{y_n\right\} \subset \mathcal{M}_\delta$ with
$
\lim _{n \rightarrow+\infty}\left|\beta_\epsilon\left(u_n\right)-y_n\right|=0 .
$
For any $u_n \in \tilde{S}(a)$, 
$$
\Upsilon_{0, T,a} \leq I_{0,T}\left(u_n\right) \leq J_{\epsilon_n,T}\left(u_n\right) \leq \Upsilon_{0,T, a}+h\left(\epsilon_n\right) \quad \forall n \in \mathbb{N},
$$
implying that
$
J_{\epsilon_n,T}\left(u_n\right) \rightarrow \Upsilon_{0, T, a} .
$
According to Lemma \ref{le4.3}, there exists $\left\{\tilde{y}_n\right\} \subset \mathbb{R}^N$ such that $y_n=\epsilon_n \tilde{y}_n $ converges to some $y \in \mathcal{M}$, and $v_n(x)=u_n\left(x+\tilde{y}_n\right)$ strongly converges to some $v \in H^s\left(\mathbb{R}^N\right)$ with $v \neq 0$. Consequently, for sufficiently large $n$, we have $\left\{y_n\right\} \subset \mathcal{M}_\delta$, and
$$
\beta_{\epsilon_n}\left(u_n\right)=y_n+\frac{\int_{\mathbb{R}^N}\left(\chi\left(\epsilon_n z+y_n\right)-y_n\right)\left|v_n\right|^2 \mathrm{d} z}{a},
$$
which leads us to conclude that
$$
\beta_{\epsilon_n}\left(u_n\right)-y_n=\frac{\int_{\mathbb{R}^N}\left(\chi\left(\epsilon_n z+y_n\right)-y_n\right)\left|v_n\right|^2 \mathrm{d} z}{a} \rightarrow 0 \quad \text { as } n \rightarrow+\infty.
$$
This completes the proof of the lemma.
\hfill{$\Box$}

{\bf Proof of Theorem 1.1.}
Decrease $\epsilon_0$ if necessary and let $\epsilon \in\left(0, \epsilon_0\right)$, then according to Lemmas \ref{le4.1}, \ref{le4.2} and \ref{le4.4}, the diagram
$$
\mathcal{M} \xrightarrow{\Phi_{\varepsilon}} \tilde{S}(a) \xrightarrow{\beta_{\varepsilon}} \mathcal{M}_\delta
$$
is well-defined. Furthermore, $\beta_\epsilon \circ \Phi_\epsilon$ is homotopic to the inclusion map id: $\mathcal{M}\rightarrow 	\mathcal{M}_\delta$. Therefore, by \cite[Lemma 4.3 ]{benci},
$
\operatorname{cat}(\tilde{S}(a)) \geq \operatorname{cat}_{\mathcal{M}_\delta}(\mathcal{M}) .
$
By the arguments presented in Lemma \ref{le2.1}, $J_{\epsilon, T}$ is bounded from below on $S(a)$. Furthermore, according to Lemma \ref{le3.5}, $J_{\epsilon, T}$ satisfies the $(P S)$ condition on $\tilde{S}(a)$.  As per the Lusternik-Schnirelmann category (see  \cite{MR1400007}), $J_{\epsilon, T}$ has at least $\operatorname{cat}(\tilde{S}(a))$ critical points on $S(a)$. Additionally, based on Lemma \ref{le3.2}, $J_{\epsilon}$ has at least  $\operatorname{cat}_{\mathcal{M}_\delta}(\mathcal{M})$ critical points on $S(a)$.


\section*{Acknowledgments}
We express our gratitude to the anonymous referee for their meticulous review of our manuscript and valuable feedback provided for its enhancement. Y. Chen was supported from the National Natural Science Foundation of China (No.12161007), Guangxi Natural Science Foundation Project (2023GXNSFAA026190), and Guangxi Science and Technology Base and Talent Project (No.AD21238019).
Z. Yang was supported by the National Natural Science Foundation of China (No.12261107), Yunnan Fundamental Research Projects (No. 202201AU070031), Scientific Research Fund of Yunnan Educational Commission (No.2023J0199) and Yunnan Key Laboratory of Modern Analytical Mathematics and Applications.

\bibliographystyle{plain}
\bibliography{nor-cho-fra}

\end{document}